# Linear and Circular Permutations with Limited Number of Repetitions

Y. Zimmels

Department of Civil and Environmental Engineering
Technion – Israel Institute of Technology
Haifa 32000, Israel

## Abstract

The problem of linear and circular permutations of $n$ identical objects in $m$ boxes, where a limit $\ell$ is imposed on the number of objects in a box, is considered. In the linear case, where the boxes are arranged as a row, two methods of solution are described. The first uses the partition diagram that is modified with prescribed number of zeros, and the second applies a direct combinatorial approach. Subject to constraints set on the relation between $n$ and $\ell$, results of the combinatorial approach are used to solve the problem where there are more than one type of objects. These solutions are then extended for the case of boxes comprising subgroups that are characterized by different values of $\ell$. In the circular case, where the boxes are arranged in the form of a circle, no direct combinatorial solution in close form is available. Consequently, in the latter case, the solution relies on the use of the modified partition diagram. This solution, which involves the application of multinomial permutations on a circle, is also extended to the case of two kinds of objects.

## Introduction

The number of permutations of *n* identical objects in *m* boxes with the number of objects, in any one of the boxes, ranging from 0 to *n*, is given by

$$I(m,n) = \binom{m+n-1}{n}. \tag{1}$$

In eq. (1), it is implicitly assumed that the option of circular equivalence, which arises when the boxes are arranged as a circle, does not exist. The latter condition is satisfied when the boxes have a linear arrangement in the form of a row. Eq. (1) also represents the case of multiple choice, or choice with repetitions. The *n* steps procedure of choice can be made from *m* kinds of sources (or populations), up to *n* times from the same source.

In this work, we consider linear and circular permutations with limited $\ell \leq n$ number of repetitions (i.e., objects in a box). We use two solution approaches; one relying on direct combinatorial analysis and the other on the Young diagram, modified with a prescribed number of zeros. Whereas both approaches lead to the same result in the linear case, only the latter seems to be readily applicable for the circular case. We do not suggest that there is no combinatorial solution for the circular case. However, such a solution (which is likely to be rather complex) still needs to be found. The derived formulas for the linear and circular permutations are extended for cases where there are more than one kind of objects and boxes that are characterized by different levels of $\ell$.

## *Statement of Problem:*

The general problem is to find the number of different ways that *n* identical objects can be distributed in *m* different boxes such that no box contains more than $\ell \leq n$ objects. The solution is sought for the linear and circular cases.

**Linear Case**

*Single type objects:*

a. Application of the partition diagram

A characteristic diagram, denoted *M*, is generated from a standard partition diagram that is modified by a prescribed number of zeros, as follows:

The partition diagram is arranged in a descending order from top to bottom and from left to right. The modified partition diagram *A*, of *n* objects (with a prescribed number of zeros so that m columns are formed) and *M* are given as,

$$
\begin{array}{ll}
n,0,\ldots\ldots,0 & m_{11}, m_{12},\ldots,m_{1t_1} \\
n-1,1,0,\ldots,0 & m_{21}, m_{22},\ldots,m_{2t_2} \\
n-2,2,0,\ldots,0 & ' \\
n-2,1,1,0,\ldots,0 & ' \hspace{3cm} (2)\\
' & ' \\
' & ' \\
' & m_{k1}, m_{k2},\ldots,m_{kt_k} \\
' & \\
1,1,1,\ldots,1,0,\ldots,0 & \\
\hspace{1cm} A & \hspace{2cm} M
\end{array}
$$

**Table 1: Diagrams *A* and *M***

Diagrams *A* and *M* each comprise *k* rows. Diagram *A* is rectangular, see Table 1. The non-negative integer $m_{ij}$ denotes the number of boxes in *A* that are occupied by the same species, i.e. number of objects. In $m_{ij}$, *i* denotes the row number, $i = 1,2,\ldots,k$ (*k* being the number of rows), and *j* the column number, or the position along the ith row, $j = 1,2,\ldots,t_i$. The number of zeros $m_{it_i}$ in the ith row is defined as,

$$m_{it_i} = m - \sum_{j=1}^{t_i-1} m_{ij}. \hspace{3cm} (3)$$

This gives a fixed sum of $m$, in each row. For example, if $m = 3$ and $n = 3$, then assuming that a box can hold up to 3 objects we have:

$$m_{11} = 1, \quad m_{12} = 2 \quad ; \quad m_{21} = 1, \quad m_{22} = 1, \quad m_{23} = 1 \quad ; \quad m_{31} = 3, \quad m_{32} = 0$$

$m_{11}$ denotes the number of boxes with three objects and $m_{12}$ the number of boxes with zeros, and so on.

The total number of permutations associated with the modified partition diagram (according to eq. (3)), is $I(m,n)$. Hence,

$$\sum_{i=1}^{k} \frac{m!}{t_i \prod_{j=1}^{} m_{ij}!} = \binom{m+n-1}{n}. \tag{4}$$

The partition diagram, which is readily constructed with standard mathematical software, is arranged in descending order of the number of objects in a box $a_{iq} (0 \leq a_{iq} \leq n, i = 1..., k, q = 1,...,m)$, in the first column and all rows: $a_{i1} \geq a_{(i+1)1}, a_{iq} \geq a_{i(q+1)}$. This facilitates the following solution to our linear case problem.

$$I(m,n,\ell) = \sum_{i=i_\ell}^{k} \frac{m!}{t_i \prod_{j=1}^{} m_{ij}!}, \tag{5}$$

where $I(m,n,\ell)$ denotes a solution for a linear case characterized by $m$ boxes, $n$ identical objects and up to $\ell \leq n$ objects in any box. The lower limit, $i = i_\ell$, denotes the first row in the partition diagram, where the number of objects in a box does not exceed $\ell$. This condition is satisfied, due to the way the diagram is structured, by all rows for which $i_\ell \leq i \leq k$. Note that in contrast to the linear case, the only available approach for solution of circular permutations (as discussed in the sequel) seems to be via the application of the modified partition diagram, as no direct combinatorial approach has been formulated to date.

b. Direct combinatorial approach

In this approach we seek the solution in the form:

$$I(m,n,\ell) = \binom{m+n-1}{n} - C, \qquad (6)$$

where $C = C(m,n,\ell) \geq 0$ denotes all permutations for which the number of objects exceeds $\ell$ at least in one box.

We start with the simplest case where $\ell < n \leq 2\ell + 1$. In this case, there can be only one box where the number of objects exceed $\ell$, i.e. from $\ell + 1$ up to $n$. Consequently, it is possible to distribute $0 \leq i \leq n - \ell - 1$ objects in the remaining $m - 1$ unoccupied boxes, while maintaining the resulting permutations in the context of $C$. As there are $m$ ways to choose the box where the number of objects will exceed $\ell$, we get

$$I(m,n,\ell) = \binom{m+n-1}{n} - m \sum_{i=0}^{n-\ell-1} \binom{m-1+i-1}{i}, \qquad (7)$$

where $\binom{m-1+i-1}{i}$ stands for the number of permutation of $i$ identical objects in $m - 1$ boxes, $i \leq \ell$. In the range $2\ell + 1 < n \leq 3\ell + 2$ the result is

$$I(m,n,\ell) = \binom{m+n-1}{n} - \left[ m \sum_{i=0}^{n-\ell-1} \binom{m-1+i-1}{i} - \binom{m}{2}\binom{m+n-2\ell-2-1}{n-2\ell-2} \right], \qquad (8)$$

where the second term in the square brackets arises from the principle of Inclusion and Exclusion[1].

Using this principle, we are now in position to generalize the result as follows,

$$I(m,n,\ell) = \binom{m+n-1}{n} - \left[ m \sum_{i=0}^{n-\ell-1} \binom{m-1+i-1}{i} - \sum_{i=2}^{j-1} (-1)^i \binom{m}{i}\binom{m+n-i\ell-i-1}{n-i\ell-i} \right], \qquad (9)$$

where $j - 1 = \lfloor n/(\ell+1) \rfloor$, which corresponds to the range $(j-1)\ell + j - 2 < n \leq j\ell + j - 1$.

Eq. (9) can be recasted as follows.

Using the identity given by Riordan[1] (p. 5, eq. 8) with redefined variables,

$$\sum_{i=0}^{n-\ell-1} \binom{m-1+i-1}{i} = \binom{m-1+n-\ell-1}{n-\ell-1} \tag{10}$$

in conjunction with eq. (9) gives,

$$I(m,n,\ell) = \binom{m+n-1}{n} - m\binom{m-1+n-\ell-1}{n-\ell-1} + \sum_{i=2}^{\lfloor n/(\ell+1) \rfloor} (-1)^i \binom{m}{i} \binom{m+n-i\ell-i-1}{n-i\ell-i}.$$

Hence,

$$I(m,n,\ell) = \sum_{i=0}^{\lfloor n/(\ell+1) \rfloor} (-1)^i \binom{m}{i} \binom{m+n-i\ell-i-1}{n-i\ell-i}. \tag{11}$$

Equation (11) has the same form given by Riordan[1], $N_n(m,\ell)$ in his notation, apart from the upper limit which he set to $m$, namely $N_n(m,\ell) = \sum_{i=0}^{m} (-1)^i \binom{m}{i} \binom{m+n-i\ell-i-1}{n-i\ell-i}$. This upper limit, which comes from his analysis of the generating function $(1+t+t^2+\ldots+t^\ell)^m$ sets in several cases we have examined $N_n(m,\ell)$ to zero, e.g. $N_6(3,2) = N_7(4,3) = N_9(4,3) = N_{11}(4,3) = 0$.

$N_6(3,2) = \binom{3}{0}\binom{3+6-1}{6} - \binom{3}{1}\binom{3+6-2-1-1}{6-2-1} + \binom{3}{2}\binom{3+6-4-2-1}{6-4-2} - \binom{3}{3}\binom{3+6-6-3-1}{6-6-3} = 0$. Derivation of eq. (11), using this generating function, is given in the Appendix. Comparing eqs. (5) and (11) gives the identity

$$\sum_{i=i_\ell}^{k} \frac{m!}{\prod_{j=1}^{t_i} m_{ij}!} = \sum_{i=0}^{\lfloor n/(\ell+1) \rfloor} (-1)^i \binom{m}{i} \binom{m+n-i\ell-i-1}{n-i\ell-i}. \tag{12}$$

By inspection of the objects arrangement in the boxes, we have $I(m, m\ell, \ell) = 1$. By virtue of eq. (11), this gives the following expansion of $\binom{m+n-1}{n}$ at $n = m\ell$,

$$\binom{m+m\ell-1}{m\ell} = 1 + \sum_{i=1}^{\lfloor n/(\ell+1) \rfloor} (-1)^{i+1} \binom{m}{i} \binom{m+m\ell-i\ell-i-1}{m\ell-i\ell-i}. \tag{13}$$

Similarly, using the fact that $I(m, m\ell - 1, \ell) = m$ gives,

$$\sum_{i=0}^{\lfloor n/(\ell+1) \rfloor} (-1)^{i} \binom{m}{i} \binom{m+m\ell-i\ell-i-2}{m\ell-i\ell-i-1} = m. \tag{14}$$

*Two kinds of objects:*

Equation (7) can be extended so as to include the case where the group of $n$ objects comprises two different subgroups, one with $n_1$ identical objects and the other with $n_2$ identical objects (different than those in the first subgroup), $n = n_1 + n_2$.

The solution for $\ell < n \leq 2\ell + 1, n_2 < n_1 \leq 2\ell + 1, 0 \leq n_2 \leq \ell, m \geq 2$, is

$$I(m, n_1, n_2, \ell) = \prod_{i=1}^{2} \binom{m+n_i-1}{n_i} - m \sum_{j=0}^{n-\ell-1} \sum_{k=0}^{\alpha} \binom{m-1+j-1}{j} \binom{m-1+k-1}{k}$$

$$= \prod_{i=1}^{2} \binom{m+n_i-1}{n_i} - m \sum_{j=0}^{n-\ell-1} \binom{m-1+j-1}{j} \binom{m-1+\alpha}{\alpha}, \tag{15}$$

where $j$ and $k$ pertain to the $n_1$ and $n_2$ subgroups, respectively, and $\alpha = \alpha(j)$ takes two different values as follows: In the range $0 \leq j \leq n - \ell - 1 - n_2$, $\alpha = n_2$, whereas $n - \ell - n_2 \leq j \leq n - \ell - 1$, sets $\alpha = n - \ell - 1 - j$. Using these definitions of $\alpha$ guarantees that $j + k \leq \ell$ is satisfied for all $j$ and $k$, so that both $j$ objects of the first kind, and $k$ of the second kind, can be distributed in the $m - 1$ boxes without constraints, as required.

Note that in derivation of the second abridged form of eq. (15), use was made of the identity (compare with eq. (10)),

$$\sum_{k=0}^{\alpha} \binom{m-1+k-1}{k} = \binom{m-1+\alpha}{\alpha}.$$

The rule, which determines the value of $\alpha$ as $j$ varies in the second form of eq. (15), remains unchanged.

Note that $n_1 > n - \ell - 1$ is imposed by the range of $n$ and $n_2 < n_1$, where eq. (15) applies. This guarantees that $j$, which denotes the number of objects (from the first subgroup) designated for distribution in the $m - 1$ boxes, can vary up to $n - \ell - 1$.

As a first example, consider the case $n_1 = \ell$, $n_2 = 1$. Here, the range $n - \ell - n_2 \le j \le n - \ell - 1$ applies, giving $n - \ell - 1 = 0$, $\alpha = 0$, and hence eq. (15) reduces to $m\binom{m+\ell-1}{\ell} - m$. These results can readily be verified by inspection.

As a second example, consider the case $n_1 = \ell$, $n_2 = 2$. This gives $n - \ell - 1 = 1, j = 0$ sets $\alpha = 1$, and $j = 1$ sets $\alpha = 0$.

Hence,

$$I(m, n_1, n_2, \ell) = \binom{m+\ell-1}{\ell}\binom{m+2-1}{2} - m\binom{m-2}{0}\left[1 + \binom{m-1}{1}\right] - m\binom{m-1}{1}$$

$$= \binom{m+1}{2}\binom{m+\ell-1}{\ell} - m(2m-1).$$

If $m = 2$, $\ell = 3$, then, $I(m = 2, n_1 = 3, n_2 = 2, \ell = 3) = 3\binom{4}{3} - 6 = 6$.

Continuing with this example, let the tuple $(n_{11}, n_{12}; n_{21}, n_{22})$ denote an allowed partition, where the 1st index in $n_{ij}$ denotes the box and the 2nd, the type of objects therein.

The following set of 6 tuples consists of the allowed partitions:

(3 , 0 | 0 , 2) , (2 , 1 | 1 , 1) , (2 , 0 | 1 , 2) , (1 , 2 | 2 , 0) , (1 , 1 | 2 , 1) , (0 , 2 | 3 , 0).

The set of 6 tuples that are not allowed is:

(3 , 2 | 0 , 0) , (3 , 1 | 0 , 1) , ( 2, 2 | 1 , 0) , ( 1 , 0 | 2 , 2) , (0 , 1 | 3 , 1) , (0 , 0 | 3 , 2).

This result agrees with the one calculated above. Equation (15) can also be extended to the case where there are $r$ groups of boxes with the ith group being characterized by $\ell_i$ (instead of the fixed $\ell$), $i = 1,...r$, and $m_i$, $\sum_{i=1}^{r} m_i = m$. Summation over all $r$ groups gives,

$$I(m,n_1,n_2,\ell,r) = \prod_{i=1}^{2} \binom{m+n_i-1}{n_i} - \sum_{i=1}^{r} m_i \sum_{j=0}^{n-\ell_i-1} \sum_{k=0}^{\alpha_i} \binom{m-1+j-1}{j}\binom{m-1+k-1}{k}$$

$$= \prod_{i=1}^{2} \binom{m+n_i-1}{n_i} - \sum_{i=1}^{r} m_i \sum_{j=0}^{n-\ell_i-1} \binom{m-1+j-1}{j}\binom{m-1+\alpha_i}{\alpha_i}, \quad (16)$$

where $0 \leq j \leq n-\ell_i-1-n_2$ sets $\alpha_i = n_2$ and $n-\ell_i-n_2 \leq j \leq n-\ell_i-1$ sets $\alpha_i = n-\ell_i-1-j$. In this case, the range where eq. (15) applies is modified as follows: $\ell_a < n \leq 2\ell_a + 1$. Here $\ell_a$ denotes the smallest value of $\ell_i, i = 1,...r$.

If $m_i = 1$ and each box has a different $\ell_i, i = 1,...,r$, then the summation turns from groups into single boxes. Observe that only in the case where $r = 1$, we have $m_i = m$. Otherwise, $m_i$ is different than $m$, which appears twice inside the double summation to the right of $m_i$. Note that if $n - \ell_i - 1 < 0$, then the corresponding term vanishes, i.e., there is no contribution from the $m_i$ boxes that are characterized by $\ell_i$, to partitions that are not allowed. This is because all $n$ objects can be placed in any one of the $m_i$ boxes without exceeding the limit $\ell_i$. If $n_2 = 0$, then $k = 0$ and eq. (16) reduces to

$$I(m,n,\ell,r) = \binom{m+n-1}{n} - \sum_{i=1}^{r} m_i \sum_{j=0}^{n-\ell_i-1} \binom{m-1+j-1}{j}$$

$$= \binom{m+n-1}{n} - \sum_{i=1}^{r} m_i \binom{m-1+n-\ell_i-1}{n-\ell_i-1}, \quad (17)$$

which represents the case of a single type objects and $r$ groups of $m_i$ boxes, each with a different limit $\ell_i, i = 1...,r$.

In the same range of $n$, eq. (15) can be further extended to cases of more than two subgroups. For example, the result for three subgroups subject to $n_1 \geq n_2 \geq n_3$, $n_1 + n_2 + n_3 = n$, is

$I(m, n_1, n_2, n_3, \ell) =$

$$\prod_{i=1}^{3} \binom{m+n_i-1}{n_i} - m \sum_{i=0}^{n-\ell-1} \sum_{j=0}^{\alpha} \sum_{k=0}^{\beta} \binom{m-1+i-1}{i}\binom{m-1+j-1}{j}\binom{m-1+k-1}{k}$$

$$= \prod_{i=1}^{3} \binom{m+n_i-1}{n_i} - m \sum_{i=0}^{n-\ell-1} \sum_{j=0}^{\alpha} \binom{m-1+i-1}{i}\binom{m-1+j-1}{j}\binom{m-1+\beta}{\beta}, \quad (18)$$

where $\alpha$ and $\beta$ take each two different values as follows: In the range $0 \leq i \leq n - \ell - 1 - n_2 - n_3$, $\alpha = n_2, \beta = n_3$, whereas in the range $n - \ell - n_2 - n_3 \leq i \leq n - \ell - 1$, $\alpha = n - \ell - 1 - i$ and $\beta = n - \ell - 1 - i - j$.

c. General solution – two kinds of objects

A general solution for the problem involving two kinds of objects is derived, using the modified partition, $A$, diagram. Let $\tilde{m}_i$ denote the number of boxes occupied by objects in the ith row of the modified $A$ diagram (eq. 2).

$$\tilde{m}_i = \sum_{j=1}^{t_i-1} m_{ij} = m - m_{it_i}, \quad (19)$$

where, as before, $m_{it_i}$ denotes the number of empty boxes in the ith row. If, as previously defined, the group of $n$ objects comprises two subgroups such that $n = n_1 + n_2$, then we denote their number of different arrangements in the $\tilde{m}_i$ boxes by $J_i$. For example, consider the case (see Table 2), $n_1 = 4, n_2 = 3, m = 8$. At i = 8 the number of different arrangements of the $n_2 = 3$ objects is $J_8 = \binom{3}{3} + \binom{2}{1} \times \binom{2}{1} + \binom{2}{1} = 7$. The term $\binom{3}{3}$ stand for

the number of ways to place 3 objects in the three boxes, one object per box. The first $\binom{2}{1}$

term represents the number of ways a pair of objects (considered as one group) can be placed in the two boxes marked in the diagram as 3,3.

The second $\binom{2}{1}$ term stands for the number of ways that the third object of $n_2$, not included in the pair, can be placed in the 2 boxes not occupied by the pair. The third $\binom{2}{1}$ term stands for the number of ways the $n_2 = 3$ objects can be placed in the first two boxes, marked 3,3 in the diagram. The second diagram, shown in Table 3, lists the different arrangements for $i = 8$ and $i = 10$. Here the number of different arrangements is 7 and 11, respectively. At $i = 10$, in agreement with the second diagram,

$J_{10} = \binom{4}{3} + \binom{2}{1} \times \binom{3}{1} + \binom{1}{1} = 11$. Using the same approach for evaluation of $J_i$ gives,

$$J_i = \binom{\tilde{m}_i}{n_2} + \binom{m(w \geq 2)}{\lfloor n_2/2 \rfloor} \times \binom{\tilde{m}_i - \lfloor n_2/2 \rfloor}{n_2 - 2\lfloor \frac{n_2}{2} \rfloor} + \binom{m(w \geq 3)}{\lfloor n_2/3 \rfloor} \times \binom{\tilde{m}_i - \lfloor n_2/3 \rfloor}{n_2 - 3\lfloor \frac{n_2}{3} \rfloor} + ..., \quad (20)$$

where $w$ denotes the number of objects in a box and m ($w \geq g$) the number of boxes with at least $g$ objects in each. The first $\binom{\tilde{m}_i}{n_2}$ term stands for the number of ways the $n_2$ objects can be distributed in the $\tilde{m}_i$ boxes, one object per box. The second $\binom{m(w \geq 2)}{\lfloor n_2/2 \rfloor}$ term represents the number of ways $\lfloor n_2/2 \rfloor$ pairs of objects can be distributed (one pair per box) in the $m$ ($w \geq 2$) boxes which carry at least 2 objects each. The third term stands for the number of different ways, the remaining object can be placed in the boxes where the pairs are absent. The following terms have similar meaning with the pairs replaced by triples and so on.

Summation over all rows $i \geq i_\ell$, for which the number of objects in a box does not exceed $\ell$, gives

$$I(m,n_1,n_2,\ell) = \sum_{i=i_\ell}^{k} J_i \frac{m!}{t_i \prod_{j=1}^{t_i} m_{ij}!}, \qquad (21)$$

where the summand in eq. (21) stands for the number of different partitions associated with the ith row. If $n_2 = 0, n_1 = n$ then $J_i = 1, i = i_\ell,...,k$ and eq. (21) reduces to eq. (5).

| $n_1 = 4$  $n_2 = 3$ , $\ell$ , $m = 8$ | |
| --- | --- |
| 7 , 0 , 0 , 0 , 0 , 0 , 0 , 0 | |
| 6 , 1 , 0 , 0 , 0 , 0 , 0 , 0 | |
| 5 , 2 , 0 , 0 , 0 , 0 , 0 , 0 | |
| 5 , 1 , 1 , 0 , 0 , 0 , 0 , 0 | |
| 4 , 3 , 0 , 0 , 0 , 0 , 0 , 0 | |
| 4 , 2 , 1 , 0 , 0 , 0 , 0 , 0 | |
| 4 , 1 , 1 , 1 , 0 , 0 , 0 , 0 | |
| ---------------------------------------- | |
| 3 , 3 , 1 , 0 , 0 , 0 , 0 , 0 | $i_\ell(\ell = 3) = 8$ |
| 3 , 2 , 2 , 0 , 0 , 0 , 0 , 0 | |
| 3 , 2 , 1 , 1 , 0 , 0 , 0 , 0 | |
| 3 , 1 , 1 , 1 , 1 , 0 , 0 , 0 | |
| ---------------------------------------- | |
| 2 , 2 , 2 , 1 , 0 , 0 , 0 , 0 | $i_\ell(\ell = 2) = 12$ |
| 2 , 2 , 1 , 1 , 1 , 0 , 0 , 0 | |
| 2 , 1 , 1 , 1 , 1 , 1 , 0 , 0 | |
| 1 , 1 , 1 , 1 , 1 , 1 , 1 , 0 | |

**Table 2: Modified Partition Diagram**

| i = 8 | | | i = 10 | | | |
|---|---|---|---|---|---|---|
| box number | | | box number | | | |
| 1 | 2 | 3 | 1 | 2 | 3 | 4 |
| 3,0 | 1,2 | 0,1 | 3,0 | 1,1 | 0,1 | 0,1 |
| 3,0 | 0,3 | 1,0 | 3,0 | 0,2 | 0,1 | 1,0 |
| 0,3 | 3,0 | 1,0 | 3,0 | 0,2 | 1,0 | 0,1 |
| 2,1 | 2,1 | 0,1 | 2,1 | 1,1 | 1,0 | 0,1 |
| 2,1 | 1,2 | 1,0 | 2,1 | 1,1 | 0,1 | 1,0 |
| 1,2 | 2,1 | 1,0 | 2,1 | 0,2 | 1,0 | 1,0 |
| 1,2 | 3,0 | 0,1 | 1,2 | 2,0 | 1,0 | 0,1 |
| | | | 1,2 | 2,0 | 0,1 | 1,0 |
| | | | 1,2 | 1,1 | 1,0 | 1,1 |
| | | | 0,3 | 2,0 | 1,0 | 1,0 |
| | | | 2,0 | 2,0 | 0,1 | 0,1 |

**Table 3: Partition of two kinds of objects in 3 and 4 boxes of rows 8 and 10, given by Table 2, respectively.**

The tuple $(u,v)$ denotes that there are $u$ objects from the first subgroup $u \leq n_1$ and $v$ objects from the second subgroup $v \leq n_2$

**Circular case**

In the circular case, the boxes are arranged in the form of a circle, keeping the same distance between adjacent boxes (so that they are uniformly distributed on the circle). In this case, we use the modified partition diagram $A$ and $M$, i.e., with prescribed number of zeros. The problem of linear multinomial permutations turns to multinomial permutations on a circle, for each row of the modified partition diagram. The number of multinomial permutations of $N$ objects on a circle, $\gamma(C_N, n^m)$, is given by[1,2]

$$\gamma\left(C_N, n^m\right) = \frac{1}{N} \sum_{d \mid \Delta} \phi(d) P\left(k^m\right) \quad, \quad P\left(k^m\right) = \frac{(k_1 + \ldots + k_m)!}{\prod_{j=1}^{m} k_j!} \quad, \quad k_j = \frac{n_j}{d} \quad, \tag{22}$$

where $x^m$, $x$ being either $n$ or $k$, denote the tuple $(x_1,\ldots,x_m)$, $\Delta = \gcd n^m$, d is any divisor of $\Delta$, $N = \sum_{i=1}^{m} n_i$, and $\phi(d)$ denotes the Euler totient function.

Application of eq. (22) to the ith row of the modified partition diagram $M$ gives the solution $\gamma_i\left(C_m, m^{t_i}\right)$ for the multinomial permutations on a circle of $m = \sum_{i=1}^{t_i} m_{ij}$ objects characterized by the tuple $m^{t_i} = \left(m_{i1},\ldots, m_{it_i}\right)$, $i = 1,2,\ldots,k$.

Thus, the counterpart of eq. (5) for the circular case is readily obtained as,

$$I_c(m,n,\ell) = \sum_{i=i_\ell}^{k} \gamma_i\left(C_m, m^{t_i}\right) , \tag{23}$$

where, as before, $i_\ell$ denotes the first row in the diagram having no more than $\ell$ objects in any of its boxes, and index $c$ denotes circularity.

Note that if $\ell = n$, than $i_n = 1$ gives the solution with no limit on the number of objects in a box. Consequently, the number of permutations, where there is at least one box with $\ell + 1$ or more objects, is given by

$$\sum_{i=1}^{i_\ell - 1} \gamma_i\left(C_m, m^{t_i}\right) = \sum_{i=1}^{k} \gamma_i\left(C_m, m^{t_i}\right) - I(m,n,\ell). \tag{24}$$

If the group of $n$ objects comprises two different subgroups, one with $n_1$ identical objects and the other with $n_2$ identical objects (different than the $n_1$), $n = n_1 + n_2$, $n_2 \leq n_1$, then, in the circular case, the counterpart solution of eq. (18) for $\ell < n \leq 2\ell + 1$, $m \geq 2$, takes the following form:

$$I_c(m,n_1,n_2,\ell) = \prod_{i=1}^{2} \sum_{i=i_\ell}^{k} \gamma\left(C_m, m^{t_i}\right) - \sum_{j=0}^{n-\ell-1} \sum_{k=0}^{\alpha} \gamma\left(C_{m-1}, m^{t_j}\right) \gamma\left(C_{m-1}, m^{t_k}\right), \quad (25)$$

where $j$ and $k$ pertain to the first and second subgroups, respectively, $\alpha = n_2$ applies for $0 \leq j \leq n-\ell-1-n_2$, whereas $n-\ell-n_2 \leq j \leq n-\ell-1$ sets $\alpha = n-\ell-1-j$. As before, this sets $j + k \leq \ell$ for all $j$ and $k$.

The second, double sum term on the right-hand-side of eq. (25), which represents the number of permutations where one box has more than $\ell$ object, vanishes at $\ell = n$ (as the upper limit of $j$ turns $-1$). In the latter case, the first term on the right-hand-side of eq. (25), which stands for the total unrestricted number of multinomial permutations on a circle of the two subgroups, gives the solution, as expected. The factor m that precedes the double sum in eq. (17) is absent from eq. (25), due to the fact that on a circle, the boxes are indistinguishable because of rotational equivalence.

The counterpart of eq. (19) for the circular case with at least two subgroups of boxes characterized by values of $\ell, \ell_1, \ell_2, ..., \ell_r, r \geq 2$,

$$I_c(m,n_1,n_2,\ell,r) = \prod_{i=1}^{2} \gamma\left(C_m, m^{t_i}\right) - \sum_{i=i}^{r} m_i \sum_{j=0}^{n-\ell_i-1} \sum_{k=0}^{\alpha_i} \gamma\left(C_{m-1}, m^{t_j}\right) \gamma\left(C_{m-1}, m^{t_k}\right) \quad (26)$$

In contrast to eq. (25), where $r = 1$ and $m_1 = m$ lead to a situation where the boxes are indistinguishable on the circle, at $r \geq 2$ (which applies to eq. (26)), the boxes are, in fact, distinguishable due to their different positions relative to other boxes on the circle. Consequently, $m_i$, which appears in the sum $\sum_{i=i}^{r}$ in eq. (26) is absent in eq. (25).

Note that the dependence of $\alpha_i = \alpha_i(j)$, on the different ranges of $j$, remains unchanged as specified for eq. (18).

The extension of eq. (25) to the case of three subgroups $(n_1, n_2, n_3)$ can be done in the same way that led to eq. (20). Note that in eq. (25), $C_{m-1}$ indicates that there are $m-1$ boxes in which the $j + k \leq \ell$ objects can be distributed, as one box must have at least $\ell + 1$ objects.

The general solution for the problem involving two kinds of objects has been derived for the linear case. Using the same arguments that led to eq. (23) gives for the circular case,

$$I_c(m, n_1, n_2, \ell) = \sum_{i=i_\ell}^{k} J_i \, \gamma\left(C_m, m^{t_i}\right). \tag{27}$$

Clearly, the same $J_i$ applies to the linear and circular cases as in both use is made of the same modified partition diagram. The difference is in the permutations within each row being $\dfrac{m!}{m_{t_i} \prod_{j=1}^{} m_{ij}}$ in the linear case and $\gamma\left(C_m, m^{t_i}\right)$ in the circular case.

# APPENDIX

Using the notation of Riordan[1], we define,

$$E(t; m, \ell) = \left(1 + t + \ldots + t^\ell\right)^m = \frac{\left(1 - t^{\ell+1}\right)^m}{(1-t)^m}. \tag{A1}$$

Expansion of the numerator and denominator gives,

$$(1-t)^{-m} = \sum_{u=0}^{\infty} \binom{m+u-1}{u} t^u, \quad \left(1 - t^{\ell+1}\right)^m = \sum_{k=0}^{m} (-1)^k \binom{m}{k} t^{k(\ell+1)}. \tag{A2}$$

Hence

$$E(t; m, \ell) = \sum_{u=0}^{\infty} \sum_{k=0}^{m} (-1)^k \binom{m}{k} \binom{m+u-1}{u} t^{u+k(\ell+1)}. \tag{A3}$$

Denoting $u = n - k(\ell + 1) \geq 0$, changes the $(k,u)$ variables to $(k,n)$, which are running in the range: $0 \leq k \leq \left\lfloor \frac{n}{\ell+1} \right\rfloor$, $0 \leq n < \infty$. Thus,

$$E(t; m, \ell) = \sum_{n=0}^{\infty} I(m, n, \ell) t^n, \tag{A4}$$

$$I(m, n, \ell) = \sum_{k=0}^{\lfloor n/(\ell+1) \rfloor} (-1)^k \binom{m}{k} \binom{m+n-k(\ell+1)-1}{n-k(\ell+1)}. \tag{A5}$$

Equations (A5) and (11) are identical except for difference in notation.